\documentclass[preprint]{elsarticle}
\usepackage{enumerate}
\usepackage{amsthm}
\usepackage{amssymb}
\usepackage{amsmath}
\usepackage{graphicx}
\usepackage{mathrsfs}
\usepackage{mathtools}
\usepackage{subfig}
\usepackage{tabularx}
\newtheorem{theorem}{Theorem}[section]
\newtheorem{prop}{Proposition}[section]
\newtheorem{cor}{Corollary}[section]
\newtheorem*{theo*}{Theorem}
\newtheorem{lm}{Lemma}[section]
\newtheorem{manualtheoreminner}{Theorem}
\newenvironment{manualtheorem}[1]{%
	\renewcommand\themanualtheoreminner{#1}%
	\manualtheoreminner
}{\endmanualtheoreminner}
\theoremstyle{definition}
\newtheorem{defin}{Definition}[section]
\newtheorem*{rem*}{Remark}
\newtheorem{rem}{Remark}[section]
\makeatletter
\def\ps@pprintTitle{%
	\let\@oddhead\@empty
	\let\@evenhead\@empty
	\def\@oddfoot{}%
	\let\@evenfoot\@oddfoot}
\makeatother
\begin{document}
	\let\today\relax
	\title{CHARACTERIZATIONS OF CONVERGENCE BY A GIVEN SET OF ANGLES IN SIMPLY CONNECTED DOMAINS}
	\author{Konstantinos Zarvalis}
	\address{Department of Mathematics, Aristotle University of Thessaloniki, 54124, Thessaloniki, Greece}
	\ead{zarkonath@math.auth.gr}
	\ead[url]{https://users.auth.gr/zarkonath}

	\begin{abstract}
		Let $\Delta$ be a simply connected domain and $f:\mathbb{D}\to\Delta$, where $\mathbb{D}$ is the unit disk, be a corresponding Riemann map. Let $\{z_n\}\subset\Delta$ be a sequence with no accumulation points inside $\Delta$. In the present article, we give necessary and sufficient conditions in terms of hyperbolic geometry which certify that $\{f^{-1}(z_n)\}$ converges to a point of $\partial\mathbb{D}$ by a certain angle $\theta$ or by a certain set of angles $[\theta_1,\theta_2]$.
	\end{abstract}
	\maketitle

	\section{Introduction}
	\par Let $\Delta$ be a simply connected domain other than $\mathbb{C}$. By the Riemann mapping theorem, this domain is conformally equivalent to the unit disk. Let $f:\mathbb{D}\to\Delta$ be a Riemann map and $\{z_n\}$ be a sequence in $\Delta$ with no accumulation points in $\Delta$. Then, $\{f^{-1}(z_n)\}$ is a sequence of the unit disk that has at least one limit point lying on $\partial\mathbb{D}$ and therefore, we can talk about the angle of the convergence. In this way, we are able to transcend the question of convergence by angle to arbitrary simply connected domains.
	\par A natural idea would be to examine the simply connected domain and its geometric properties to provide an answer about the angle of the convergence of the preimage sequence in the unit disc. However, the domain $\Delta$ can be extremely complicated, thus making such a procedure very difficult and at times impossible.
	\par In \cite{orth}, Bracci, Contreras, D\'{i}az-Madrigal and Gaussier introduced an interesting technique to tackle the above problem. Instead of looking at the given simply connected domain $\Delta$, we find another one that contains $\Delta$, has better geometric properties and is easier to work with. More precisely, we find a simply connected domain $U$ such that $\Delta\subset U$ and $\Delta$ contains a suitable horodisk of $U$.
	\par To be more concrete, we explain first the notion of a horodisk. Let $\sigma\in\partial\mathbb{D}$ and $R>0$. Then, the \textit{horodisk of} $\mathbb{D}$ \textit{of radius} $R$ \textit{centered at} $\sigma$ is the set $$E_\mathbb{D}(\sigma,R)=\{z\in\mathbb{D}:|\sigma-z|^2<R(1-|z|^2)\}.$$ Again, we need to generalize this definition. Let $U\subsetneq\mathbb{C}$ be a simply connected domain, $\partial_C U$ the set of its prime ends and $f:\mathbb{D}\to U$ a Riemann map. Suppose that $\xi\in\partial_C U$ and that the prime end $\xi$ corresponds through $f$ to the point $\sigma\in\partial\mathbb{D}$. Then, the \textit{horodisk of $U$ of radius $R$ centered at $\xi$} is the set $$E_U(\xi,R)=f(E_\mathbb{D}(\sigma,R)).$$
	\par The ``suitable horodisk'' is one that is centered at the limit point of the sequence. Setting $\xi$ the limit point and using the above terminology, we have $E_U(\xi,R)\subset\Delta\subseteq U$ (the equality is allowed for the case when $\Delta$ is already simple). Therefore, the given simply connected domain $\Delta$ is restricted between $U$ and its horodisk and as the sequence converges to $\xi$ the geometry of $\Delta$ seems more and more similar to that of the simpler domain $U$.
	\par So, the central idea is the use of a domain with simpler geometry and its horodisks in order to find a necessary and sufficient condition for convergence by a certain angle. In \cite{orth}, the authors deal only with the case of \textit{orthogonal} convergence, when the angle of convergence is $\frac{\pi}{2}$, and provide the corresponding theorem. In the present article, we will utilise their method to produce an analogue result for any angle-set $[\theta_1,\theta_2]\subset(0,\pi)$.
	\par Moreover, when we speak of the geometry of these simply connected domains, we mostly talk in terms of hyperbolic geometry. So, before mentioning the theorem in \cite{orth}, we first give some useful definitions to render its statement clearer. To begin with, the \textit{hyperbolic distance} in $\mathbb{D}$ is given by (see e.g. \cite[Section 1.3]{book}) $$k_\mathbb{D}(z,w)=\frac{1}{2}\log\frac{1+\left|\frac{z-w}{1-\bar{z}w}\right|}{1-\left|\frac{z-w}{1-\bar{z}w}\right|},\;\;z,w\in\mathbb{D}.$$ Now, let $\Delta\subsetneq\mathbb{C}$ be a simply connected domain and $f:\Delta\to\mathbb{D}$ a corresponding Riemann map. We denote by $k_\Delta$ the hyperbolic distance in $\Delta$ given by $$k_\Delta(z,w)=k_\mathbb{D}(f(z),f(w)),\;\;z,w\in\Delta.$$ It can be verified that this definition is indeed independent of the choice of the conformal map $f$. Moreover, if $T\subset\mathbb{R}$ is a segment (open, closed, semi-open) and $\gamma:T\to\Delta$ is a smooth curve, then we denote $$k_\Delta(z,\gamma)=\inf\limits_{t\in T}k_\Delta(z,\gamma(t)),\;z\in\Delta.$$
	\par Now, we are ready for the theorem.
	\begin{theo*}[\cite{orth}]
		Let $\Delta\subsetneq\mathbb{C}$ be a simply connected domain and $f:\mathbb{D}\to\Delta$ a Riemann map. Let $\{z_n\}\subset\Delta$ be a sequence with no accumulation points in $\Delta$. Then, there exists a $\sigma\in\partial\mathbb{D}$ such that $\{f^{-1}(z_n)\}$ converges orthogonally to $\sigma$ if and only if there exist a simply connected domain $U\subsetneq\mathbb{C},\; \xi\in\partial_C U$ and $R>0$ such that
		\begin{enumerate}[\normalfont(i)]
			\item $E_U(\xi,R)\subset\Delta\subseteq U$,
			\item $\lim\limits_{n\to+\infty}k_U(z_n,\gamma)=0$, where $\gamma:[0,+\infty)\to U$ is any geodesic for the hyperbolic distance in $U$ such that $\lim\limits_{t\to+\infty}\gamma(t)=\xi$ in the Carath\'{e}odory topology of $U$.
		\end{enumerate} 
	\end{theo*}
	\par We will generalize this theorem to give a result for convergence by any angle $\theta\in(0,\pi)$, by providing a corollary to our main theorem which will be an even more general result concerning the convergence by any angle-set $[\theta_1,\theta_2]\subset(0,\pi)$.

	\par In order to obtain the proof, we will mostly utilize techniques of harmonic measure; techniques upon which we are going to elaborate further in the next section. The corresponding proof given in \cite{orth} for orthogonal convergence relies on properties of the hyperbolic distance. 
	\par We need to provide a rigid definition of convergence by angle.
	\begin{defin}
		A sequence $\{z_n\}\subset\mathbb{D}$ is said to converge \textit{by angle $\theta$}, $\theta\in[0,\pi],$ to a point $\sigma\in\partial\mathbb{D}$, provided $\lim\limits_{n\to+\infty}z_n=\sigma$ and $\lim\limits_{n\to+\infty}\arg(1-\bar{\sigma}z_n)=\frac{\pi}{2}-\theta$. In a similar manner, a continuous curve $\gamma:[0,+\infty)\to\mathbb{D}$ converges \textit{by angle $\theta$}, $\theta\in[0,\pi]$, to a point $\sigma\in\partial\mathbb{D}$, provided $\lim\limits_{t\to+\infty}\gamma(t)=\sigma$ and $\lim\limits_{t\to+\infty}\arg(1-\bar{\sigma}\gamma(t))=\frac{\pi}{2}-\theta$.
	\end{defin}
	\par In the extreme cases when $\theta=0$ or $\theta=\pi$, the convergence is characterized as \textit{tangential}. Nevertheless, a sequence in the unit disk may oscillate while converging to its limit on the unit circle and as a result this sequence can have various (and even infinite) subsequences that converge to the limit by different angles. In this case, we cannot say that the sequence converges by a certain angle. We are in need of a new, more inclusive definition.
	
	\begin{defin}
		A sequence $\{z_n\}\subset\mathbb{D}$ is said to converge \textit{by angle-set} $[\theta_1,\theta_2]\subset[0,\pi]$ to a point $\sigma\in\partial\mathbb{D}$, provided $\lim\limits_{n\to+\infty}z_n=\sigma$ and the cluster set of $\arg(1-\bar{\sigma}z_n)$ is $\left[\frac{\pi}{2}-\theta_2,\frac{\pi}{2}-\theta_1\right]$.
	\end{defin} 
	
	\par In this article, we will only entertain the possibility when $[\theta_1,\theta_2]\subset(0,\pi)$. Clearly, as was the case previously,  the definition can be extended to a continuous curve $\gamma:[0,+\infty)\to\mathbb{D}$.
	\par It is a well known fact that the cluster set of $\arg(1-\bar{\sigma}\gamma(t))$ for a continuous curve $\gamma$ is a compact, connected subset of $[0,\pi]$. So, this cluster set must be either a singleton $\{\theta\}$ or a continuum $[\theta_1,\theta_2]$. Consequently, a continuous curve falls into one of the two categories described by Definitions 1.1 and 1.2. However, for a sequence $\{z_n\}$, things can get much more complicated. For example, the cluster set of $\arg(1-\bar{\sigma}z_n)$ could be a finite set of angles, or even a sequence of angles. In particular, if the cluster set of $\arg(1-\bar{\sigma}z_n)$ is a subset of $(0,\pi)$, the convergence is called \textit{non-tangential}. For the purposes of this article, we only examine the first two categories.
	\par The goal of our article is to provide a necessary and sufficient condition for convergence by angle-set. To this direction, we first need the notion of hyperbolic sectors. Let $\Delta$  be a simply connected domain and $\gamma:[\alpha,+\infty)\to\Delta,\;\alpha\ge-\infty$, a geodesic for the hyperbolic distance in $\Delta$ with the added property $\lim\limits_{t\to+\infty}k_\Delta(\gamma(t),\gamma(t_0))=+\infty$, for some $t_o\in(a,+\infty)$. For $R>0$, we call the set $$S_\Delta(\gamma,R)=\{z\in\Delta:k_\Delta(z,\gamma)<R\}$$ a \textit{hyperbolic sector around} $\gamma$ \textit{of amplitude} $R$.
	\par In \cite{non-tan}, again Bracci, Contreras, D\'{i}az-Madrigal and Gaussier proved the following theorem that characterizes non-tangential convergence:
	
	\begin{theo*}[\cite{non-tan}]
		Let $\Delta\subsetneq\mathbb{C}$ be a simply connected domain and let $f:\mathbb{D}\to\Delta$ be a Riemann map. Let $\{z_n\}\subset\Delta$ be a sequence with no accumulation points in $\Delta$. Then, $\{f^{-1}(z_n)\}$ converges non-tangentially to a point $\sigma\in\partial\mathbb{D}$ if and only if there exist a simply connected domain $U\subsetneq\mathbb{C}$, a geodesic $\gamma:[0,+\infty)\to U$ of $U$ such that $\lim\limits_{t\to+\infty}k_U(\gamma(t),\gamma(0))=+\infty$ and $R>R_0>0$ such that
		\begin{enumerate}[\normalfont(i)]
			\item $S_U(\gamma,R)\subset\Delta\subseteq U$,
			\item there exists $n_0\ge0$ such that $z_n\in S_U(\gamma,R_0)$ for all $n\ge n_0$.
		\end{enumerate}
	\end{theo*}
	
	\par From their theorem, we borrow the central idea of our characterization of convergence by angle-set. This idea is the use of hyperbolic sectors and we combine this idea again with the use of a simple, in terms of geometry, simply connected domain and its horodisks. However, the statement of the theorem demands certain new definitions and notations which will be dealt with in detail in Section 3. 
	\par Briefly, let $\{z_n\}\subset\Delta$ be a sequence with no accumulation points in $\Delta$, where $\Delta$ is a simply connected domain, and $\lim\limits_{n\to+\infty}z_n=\xi\in\partial_C\Delta$. We will introduce the sets $A_\Delta(\gamma,\theta_1,\theta_2)$, where $\gamma:[0,+\infty)\to\Delta$ is a geodesic for the hyperbolic distance in $\Delta$ with $\lim\limits_{t\to+\infty}\gamma(t)=\xi$ and $0<\theta_1<\theta_2<\pi$ are the angles that define the angle-set $[\theta_1,\theta_2]$ of convergence. The definition of these sets relies exclusively on the geodesic $\gamma$ and its hyperbolic sectors with amplitudes depending on the angles. We will say that the sequence $\{z_n\}$ \textit{exhausts} the set $A_\Delta(\gamma,\theta_1,\theta_2)$ if it is eventually contained in every set $A_\Delta(\gamma,\omega_1,\omega_2),\; \omega_1<\theta_1,\omega_2>\theta_2$ and, in a way that will be explained later, eventually ``fills'' the set $A_\Delta(\gamma,\theta_1,\theta_2)$. Our theorem is the following.
	
	\begin{theorem}
		Let $\Delta\subsetneq\mathbb{C}$ be a simply connected domain and $f:\mathbb{D}\to\Delta$ be a Riemann map. Let $\{z_n\}\subset\Delta$ be a sequence with no accumulation points in $\Delta$. Then, there exists a $\sigma\in\partial\mathbb{D}$ such that $\{f^{-1}(z_n)\}$ converges by angle-set $[\theta_1,\theta_2]\subset(0,\pi)$ to $\sigma$ if and only if there exist a simply connected domain $U\subsetneq\mathbb{C}$, $\xi\in\partial_C U$, $R>0$ and a geodesic $\gamma:[0,+\infty)\to U$ for the hyperbolic distance in $U$ with $\lim\limits_{t\to+\infty}\gamma(t)=\xi$ in the Carath\'{e}odory topology of $U$ such that
		\begin{enumerate}[\normalfont(i)]
			\item $E_U(\xi,R)\subset\Delta\subseteq U,$
			\item $\gamma(0)\in\Delta,$
			\item $\{z_n\}$ exhausts the set $A_U(\gamma,\theta_1,\theta_2).$
		\end{enumerate}
	\end{theorem}

	\par Again, our proof relies heavily on harmonic measure and the conformally invariant nature of hyperbolic distance and harmonic measure. In \cite{non-tan}, the authors rely on \textit{Gromov's hyperbolicity theory} and \textit{quasi-geodesics}. Our proof thus shows that Gromov's theory is not necessary. However, Gromov's theory and the notion of quasi-geodesics are certainly very useful for a profound understanding of the problem.
	
	\par Finally, in Section 4, we will give an application on continuous semigroups of holomorphic functions of the unit disk.

\section{Preliminaries}
\subsection{Hyperbolic Distance}
\par Before moving on to the main body of the article, we must first review some basic facts and properties concerning two conformal invariants: hyperbolic distance and harmonic measure.
\par Starting with the hyperbolic distance, its property of conformal invariance is directly implied by its definition with respect to an arbitrary simply connected domain, other than the complex plane. In addition, hyperbolic distance has a monotonicity property: let $\Delta_1,\Delta_2$ be two simply connected domains satisfying $\Delta_1\subset\Delta_2$. Then, for all $z,w\in\Delta_1$ we have that $$k_{\Delta_2}(z,w)\le k_{\Delta_1}(z,w).$$

\subsection{Harmonic Measure}
\par Next, we turn to \textit{harmonic measure}, which is conformally invariant as well. A comprehensive presentation of its theory can be found in \cite{rans}. Let $\Delta$ be a domain in $\mathbb{C}$ with non-polar boundary. Let $B$ be a Borel subset of $\partial\Delta$. Then, the harmonic measure of $B$ with respect to $\Delta$ is exactly the solution of the generalized Dirichlet problem for the Laplacian in $\Delta$ with boundary function equal to $1$ on $B$ and to $0$ on $\partial\Delta\setminus B$. For the harmonic measure of $B$ with respect to $\Delta$ and for $z\in\Delta$ we use the notation $\omega(z,B,\Delta)$. It is known that for a fixed $z\in\Delta$, $\omega(z,\cdot,\Delta)$ is a Borel probability measure on $\partial\Delta$.
\par A useful property of harmonic measure is its domain monotonicity. In detail, let $\Delta_1\subset\Delta_2,\;B\subset\partial\Delta_1\cap\partial\Delta_2$ Borel. Then, for all $z\in\Delta_1$, we have $$\omega(z,B,\Delta_1)\le\omega(z,B,\Delta_2).$$ The above statement can be made more precise by means of what is called the \textit{strong Markov property} of harmonic measure (see e.g. \cite[p.307]{conway}). Again, let $\Delta_1\subset\Delta_2$ and $B\subset\partial\Delta_1\cap\partial\Delta_2$ Borel. Then, for $z\in\Delta_1,$ $$\omega(z,B,\Delta_2)=\omega(z,B,\Delta_1)+\int\limits_{\partial\Delta_1\setminus\partial\Delta_2}\omega(\zeta,B,\Delta_2)\cdot\omega(z,d\zeta,\Delta_1).$$
\begin{rem}
	 Let $B$ be an arc of $\partial\mathbb{D}$ with endpoints $a,b$. Then, we know (see e.g. \cite[p.155]{cara}) that the level set $$L_k=\{z\in\mathbb{D}:\omega(z,B,\mathbb{D})=k\},\;0<k<1,$$ is a circular arc (or a diameter in case $B$ is a half-circle and $k=\frac{1}{2}$) in $\mathbb{D}$ with endpoints $a,b$ that meets $\partial\mathbb{D}$ with angle $k\pi$. Therefore, a sequence $\{z_n\}\subset\mathbb{D}$ that converges to $a$ or $b$, converges by angle $\theta$ if and only if $\lim\limits_{n\to+\infty}\omega(z_n,B,\mathbb{D})=\frac{\theta}{\pi}$. The analogous result is also true for any curve converging to $a$ or $b$ as well. Certainly, inside the subdomain of the unit disk bounded by the level set $L_k$ and the arc $B$, the harmonic measure $\omega(z,B,\mathbb{D})$ is larger than $k$, while inside the complementary subdomain, the corresponding harmonic measure is smaller than $k$.
\end{rem}

\par Finally, we refer to \cite{rans} for some results about harmonic measure on certain particular domains. Denoting by $U_1$ the upper half plane $\{z:\Im z>0\}$, we have $$\omega(z,[a,b],U_1)=\frac{1}{\pi}\arg\left(\frac{z-b}{z-a}\right),$$ where $[a,b]\subset\mathbb{R}, z\in U_1$. Then, for any sequence $\{z_n\}\subset U_1$ with $\lim\limits_{n\to+\infty}z_n=\infty$, it is easily calculated by the above formula that $\lim\limits_{n\to+\infty}\omega(z_n,[a,b],U_1)=0$. Furthermore, denoting by $U_2$ the set $\{z:a<\arg z<\beta\}$, where $(a,\beta)\subset(-\pi,\pi)$, for $z\in U_2$ we get $$\omega(z,\{\arg z=\beta\},U_2)=\frac{\arg z-a}{\beta-a}.$$ A direct consequence of the last relation is that for the right half plane $\mathbb{H}=\{z:\Re z>0\}$ we have $$\omega(t,\{\arg z=-\frac{\pi}{2}\},\mathbb{H})=\omega(t,\{\arg z=\frac{\pi}{2}\},\mathbb{H})=\frac{1}{2},$$ for all positive $t$. Of course, using conformal mappings and combining the above equalities, we can find numerous other results.

\subsection{Convergence by Angle}
\par It is necessary to understand the elementary origin of Definitions 1.1 and 1.2 and comprehend  what they signify. To this goal, let $\sigma\in\partial\mathbb{D}$ and $T$ be the tangent of $\partial\mathbb{D}$ at $\sigma$. Consider $\Gamma$ to be a straight line that intersects $\partial\mathbb{D}$ at $\sigma$. Restricting ourselves to the half-plane determined by $T$ and $\mathbb{D}$, the lines $T$ and $\Gamma$ form two angles. We call $\theta$ the angle formed by $\Gamma$ and the ray of $T$ extending to the right of $\sigma$ (where right and left are determinded by considering the half-plane we chose as the lower one, as shown in Figure 1). Let $\{z_n\}\subset\Gamma\cap\mathbb{D}$ with $\lim\limits_{n\to+\infty}z_n=\sigma$. Geometrically speaking, we could say that $\{z_n\}$ converges by angle $\theta$ to $\sigma$. To be more rigid, for all $n\in\mathbb{N}$, we can write $z_n=\sigma+r_ne^{i\theta_n},\;r_n>0,\;\theta_n\in[0,2\pi)$. Looking at Figure 1 and after doing some elementary computations, it is easily verified that for every term of the sequence $\{z_n\}$, we have $$\arg(1-\bar{\sigma}z_n)=\frac{\pi}{2}-\theta,$$ 
and of course the same is true for the limit as $n\to+\infty$. As a result, the practical meaning of Definition 1.1 is that for a sequence $\{z_n\}$ to converge by angle $\theta$ to $\sigma$, it is necessary that there is a $n_0\in\mathbb{N}$ such that for all $n>n_0$, all the terms $z_n$ are sufficiently close to the line $\Gamma$. A similar conclusion holds in the case of a continuous curve.

\begin{figure}
	\centering
	\includegraphics[scale=0.4]{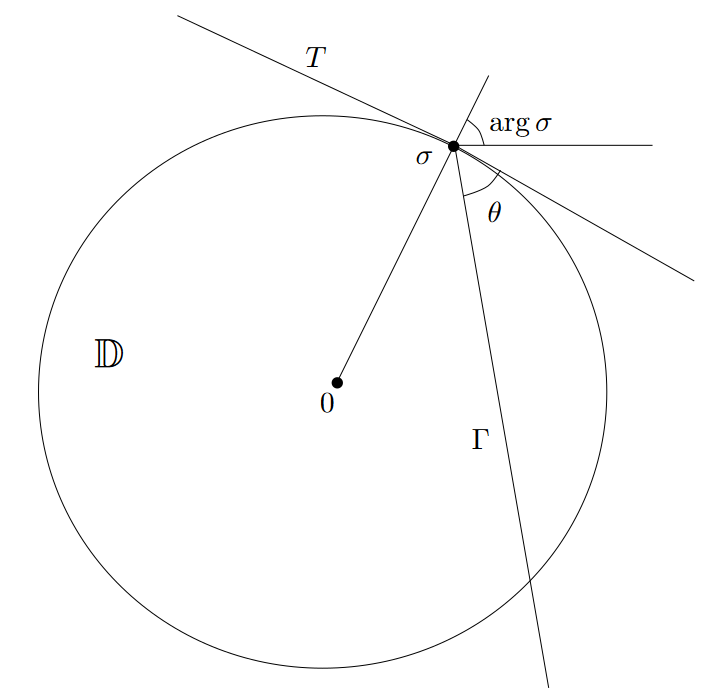}
	\caption{Origin of Definition 1.1}
\end{figure}

\section{Convergence by angle-set}
\par Being restricted to $(0,\pi)$, the convergence by angle-set is non-tangential. Therefore, if a sequence $\{z_n\}\subset\mathbb{D}$ converges to a point $\sigma\in\partial\mathbb{D}$ by a specific angle-set, it is necessary that this sequence is eventually contained in a \textit{Stolz region of vertex} $\sigma$, that is the set $$\{z\in\mathbb{D}:|\sigma-z|<R(1-|z|)\},$$ for some $R>1$. Of course, the inverse implication is also true. It is known that if we take a geodesic of $\mathbb{D}$ with endpoint at  $\sigma$, then every Stolz region of vertex $\sigma$ is equivalent to a hyperbolic sector around $\gamma$, a notion we have already introduced. Thus, we could roughly say that convergence by angle-set might be characterized with the help of hyperbolic sectors. Indeed, the principal role of this section will be played by hyperbolic sectors. We will, first, need some lemmas and definitions.

\par We commence with a known fact about geodesics. Let $\Delta\subsetneq\mathbb{C}$ be a simply connected domain and $\gamma:(a,b)\to\Delta$, where $a,b$ can be $-\infty,\infty$ respectively, a geodesic for the hyperbolic distance in $\Delta$. We suppose that $\lim\limits_{t\to a^+}\gamma(t),\lim\limits_{t\to b^-}\gamma(t)\in\partial\Delta$ or, in other words, that with this parametrization we have the whole geodesic and not only part of it. Then, the set $\Delta\setminus\gamma((a,b))$ consists of two simply connected components.

\par Using the above remark, the following notation can be stated naturally. Let $\gamma:(-1,1)\to\mathbb{D}$ be a geodesic in $\mathbb{D}$ such that $\lim\limits_{t\to -1^+}\gamma(t)=\tau\in\partial\mathbb{D}$ and $\lim\limits_{t\to 1^-}\gamma(t)=\sigma\in\partial\mathbb{D}$. Then, $\mathbb{D}\setminus\gamma((-1,1))$ consists of two simply connected components. Consider $\gamma$ to be oriented with direction towards $\sigma$. We denote by $\mathbb{D}_\sigma^+$ the simply connected component of $\mathbb{D}\setminus\gamma((-1,1))$ that lies to the right of $\gamma$. We denote by $\mathbb{D}_\sigma^-$ the one that lies to the left.

\par Now, we need some more important notations to render everything that will follow more legible. Let $[\theta_1,\theta_2]\subset(0,\pi)$, $\sigma\in\partial\mathbb{D}$ and $\gamma:[0,1)\to\mathbb{D}$ be a geodesic (actually part of it) in $\mathbb{D}$ with $\lim\limits_{t\to 1^-}\gamma(t)=\sigma$. Of course, the geodesic $\gamma$ can be extended in a unique way in order to reach its second endpoint on $\partial\mathbb{D}$. Through this extension, we can naturally obtain, again, the sets $\mathbb{D}_\sigma^-$ and $\mathbb{D}_\sigma^+$. 
From now on, we set
$$R(\theta)=k_\mathbb{H}(1,e^{i\theta}),\;\theta\in(-\frac{\pi}{2},\frac{\pi}{2}).$$
The importance of this number will be understood better later on. For now, we can mention some properties. First of all, it is obvious that $R(0)=0$. In addition, because of the symmetry of the right half-plane with respect to the real axis, it is easily seen that $R(\theta)=R(-\theta)$. Moreover, vertical translations and scalings are conformal automorphisms of $\mathbb{H}$. Therefore, $$R(\theta)=k_\mathbb{H}(1+iy,e^{i\theta}+iy)=k_\mathbb{H}(r,re^{i\theta}),$$ for all $y\in\mathbb{R}$ and for all $r>0$. Composing a vertical translation and a scaling we can get new combinations. Finally, using the Cayley transform $C:\mathbb{H}\to\mathbb{D}$, with $C(z)=\frac{z-1}{z+1}$, we get $$R(\theta)=k_\mathbb{D}\left(0,\frac{e^{i\theta}-1}{e^{i\theta}+1}\right)=k_\mathbb{D}\left(0,\left|\frac{e^{i\theta}-1}{e^{i\theta}+1}\right|\right).$$ Using the formula for the hyperbolic distance in the unit disk, we can find that $R(\theta)=\textrm{arctanh}(|\tan (\frac{\theta}{2})|)$.
\par Next, we mention the following lemma.

\begin{lm}[\cite{non-tan}]
	Let $\gamma:[0,+\infty)\to\mathbb{H}$ be a geodesic for the hyperbolic distance in $\mathbb{H}$ such that $\gamma([0,+\infty))=[r_0,+\infty),\;\gamma(0)=r_0$, for some $r_0>0$. Then, for every $R>0$, there exists $\beta\in(0,\frac{\pi}{2})$, with $k_\mathbb{H}(r_0,r_0e^{i\beta})=R$, such that $$S_\mathbb{H}(\gamma,R)=\{\rho e^{i\theta}:\rho>r_0,|\theta|<\beta\}\cup\{z\in\mathbb{H}:k_\mathbb{H}(r_0,z)<R\}.$$
\end{lm}

\begin{figure}
	\centering
	\includegraphics[scale=0.7]{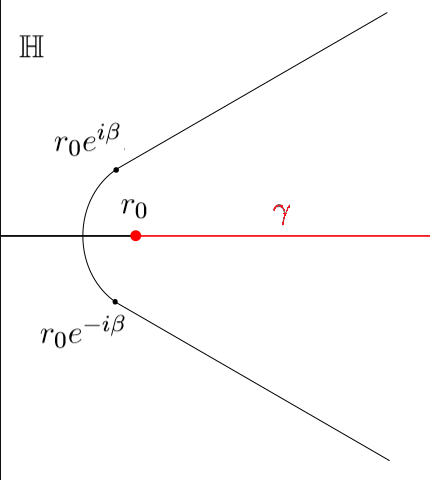}
	\caption{Hyperbolic sector in the right half-plane}
\end{figure}

\par With the help of Lemma 3.1, we can see (Figure 2) that a hyperbolic sector in the right half-plane around a ray on the positive semi-axis is actually a simply connected domain bounded by a circular arc that joins the two points $r_0e^{i\beta}, r_0e^{-i\beta}$ and by the two rays $\{\rho e^{i\beta}:\rho\ge r_0\},\{\rho e^{-i\beta}:\rho\ge r_0\}$. Also, it is symmetric with respect to the real line. Thanks to conformal invariance, the right half-plane and the ease in picturing its hyperbolic sectors will play a crucial role later on.

\par Now, we can picture a hyperbolic sector of amplitude $R(\theta)$. Let $\gamma:[0,+\infty)\to\mathbb{H}$ be a geodesic of the right half-plane such that its image is a ray on the positive real axis with $\gamma(0)=r_0$. Then, utilising the above lemma, for $\theta_0\in(-\frac{\pi}{2},\frac{\pi}{2})\setminus\{0\}$ we have 
$$S_\mathbb{H}(\gamma,R(\theta_0))=\{\rho e^{i\theta}:\rho>r_0,\;|\theta|<|\theta_0|\}\cup\{z\in\mathbb{H}:k_\mathbb{H}(r_0,z)<R(\theta_0)\}.$$
So, the sector $S_\mathbb{H}(\gamma,R(\theta_0))$ eventually contains every ray $\{z:\arg z=\theta\}$, for $\theta\in(-|\theta_0|,|\theta_0|)$. Using harmonic measure, this means that the sector $S_\mathbb{H}(\gamma,R(\theta_0))$ eventually contains every level set $\{z:\omega(z,\{\arg z=\frac{\pi}{2}\},\mathbb{H})=\frac{\theta+\frac{\pi}{2}}{\pi}\}$, for $\theta\in(-|\theta_0|,|\theta_0|)$. But this gives us some useful information for convergence by angle. Through the Cayley transform that corresponds $\infty$ to $1$, all the above certify that this sector eventually contains sequences whose images converge to $1$ by angles inside $(-|\theta_0|+\frac{\pi}{2},|\theta_0|+\frac{\pi}{2})$.
\par Moreover, for any $y\in\mathbb{R}$ and for $\theta\in(-\frac{\pi}{2},\frac{\pi}{2})$ fixed, we can easily compute that 
$$\lim_{\substack{\zeta\to\infty \\ \arg\zeta=\theta}}\omega(\zeta+iy,\{\arg z=\frac{\pi}{2}\},\mathbb{H})=\frac{\theta+\frac{\pi}{2}}{\pi}.$$ Of course, the ray $t+iy,\;t\ge r_0$, where $r_0>0$, is a geodesic of $\mathbb{H}$. Considering $\gamma(t)=r_0+t+iy$, then the hyperbolic sector $S_\mathbb{H}(\gamma,R(\theta_0))$ again eventually contains all the rays $\{z:\arg z=\theta\}+iy$, for $\theta\in(-|\theta_0|,|\theta_0|)$. Therefore, for any $\theta\in(-|\theta_0|+\frac{\pi}{2},|\theta_0|+\frac{\pi}{2})$, any sequence such that its image in $\mathbb{D}$ converges to $1$ by angle $\theta$ is contained in some sector $S_\mathbb{H}(\gamma,R(\theta_0))$, where $\gamma:[0,+\infty)\to\mathbb{H}$ is a geodesic of $\mathbb{H}$ such that $\lim\limits_{t\to+\infty}\gamma(t)=\infty$. Through the conformal invariance of the hyperbolic distance, a similar thing can be stated for any simply connected domain.
\par Finally, combining all the above information, we can even predict the following: any sequence of $\mathbb{H}$ such that its image converges to $1$ by angle-set $[\theta_1,\theta_2]$ must be eventually contained in a sector $S_\mathbb{H}(\gamma,R(\frac{\pi}{2}-\theta_1))$ and not be eventually contained in the sector $S_\mathbb{H}(\gamma,R(\frac{\pi}{2}-\theta_2))$ or the inverse, depending on the angles.

So, we are in a position to define the set $A_\mathbb{D}(\gamma,\theta_1,\theta_2)$, for $[\theta_1,\theta_2]\subset(0,\pi)$, where $\gamma:[0,1)\to\mathbb{D}$ is a geodesic in $\mathbb{D}$ such that $\lim\limits_{t\to 1^-}\gamma(t)=\sigma$. But first, we ought to define the auxiliary sets $\tilde{A}_\mathbb{D}(\gamma,\theta_1,\theta_2)$ as follows:

	\begin{enumerate}[(i)]
		\item if $[\theta_1,\theta_2]\subset(0,\frac{\pi}{2})$, then \\
		$\tilde{A}_\mathbb{D}(\gamma,\theta_1,\theta_2)=\left[S_\mathbb{D}(\gamma,R(\frac{\pi}{2}-\theta_1))\setminus S_\mathbb{D}(\gamma,R(\frac{\pi}{2}-\theta_2))\right]\cap\mathbb{D}_\sigma^+$,
		\item if $[\theta_1,\theta_2]\subset(\frac{\pi}{2},\pi)$, then \\
		$\tilde{A}_\mathbb{D}(\gamma,\theta_1,\theta_2)=\left[S_\mathbb{D}(\gamma,R(\frac{\pi}{2}-\theta_2))\setminus S_\mathbb{D}(\gamma,R(\frac{\pi}{2}-\theta_1))\right]\cap\mathbb{D}_\sigma^-$,
		\item if $[\theta_1,\theta_2]=[\theta,\pi-\theta]$, where $\theta\in(0,\frac{\pi}{2})$, then \\
		$\tilde{A}_\mathbb{D}(\gamma,\theta_1,\theta_2)=S_\mathbb{D}(\gamma,R(\frac{\pi}{2}-\theta_1))= S_\mathbb{D}(\gamma,R(\frac{\pi}{2}-\theta_2))$,
		\item if $[\theta_1,\theta_2]=[\theta_1,\frac{\pi}{2}]$, then \\
		$\tilde{A}_\mathbb{D}(\gamma,\theta_1,\theta_2)=S_\mathbb{D}(\gamma,R(\frac{\pi}{2}-\theta_1))\cap\mathbb{D}_\sigma^+$,
		\item if $[\theta_1,\theta_2]=[\frac{\pi}{2},\theta_2]$, then \\
		$\tilde{A}_\mathbb{D}(\gamma,\theta_1,\theta_2)=S_\mathbb{D}(\gamma,R(\frac{\pi}{2}-\theta_2))\cap\mathbb{D}_\sigma^-$,
		\item if $[\theta_1,\theta_2]$ satisfies $\theta_1<\frac{\pi}{2}<\theta_2$, then \\
		$\tilde{A}_\mathbb{D}(\gamma,\theta_1,\theta_2)=\left[S_\mathbb{D}(\gamma,R(\frac{\pi}{2}-\theta_1))\cap\mathbb{D}_\sigma^+\right]\cup\left[S_\mathbb{D}(\gamma,R(\frac{\pi}{2}-\theta_2))\cap\mathbb{D}_\sigma^-\right]$.
	\end{enumerate}

Finally, for any of the above six cases, we define
$$A_\mathbb{D}(\gamma,\theta_1,\theta_2):=\overline{\tilde{A}_\mathbb{D}(\gamma,\theta_1,\theta_2)}\cap\mathbb{D},$$
where the overline signifies the closure of the set. In other words, the set $A_\mathbb{D}(\gamma,\theta_1,\theta_2)$ is just the closure of the auxiliary set $\tilde{A}_\mathbb{D}(\gamma,\theta_1,\theta_2)$ minus the point $\sigma$. Obviously, the definition of the sets $\tilde{A}_\mathbb{D}(\gamma,\theta_1,\theta_2)$ is not really needed, but we try to avoid the repetitive use of the closure in all of the six cases.

\begin{figure}
	\centering
	\def\tabularxcolumn#1{m{#1}}
	\begin{tabularx}{\linewidth}{@{}cXX@{}}
		\begin{tabular}{cc}
			\subfloat[\({[\theta_1,\theta_2]}\subset(0,\frac{\pi}{2})\)]{\includegraphics[scale=0.3]{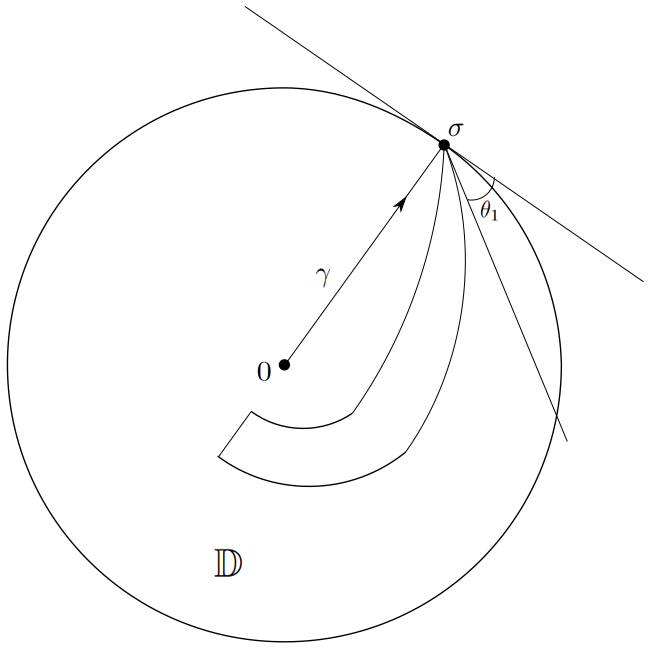}}
			& \subfloat[\({[\theta_1,\theta_2]}\subset(\frac{\pi}{2},\pi)\)]{\includegraphics[scale=0.36]{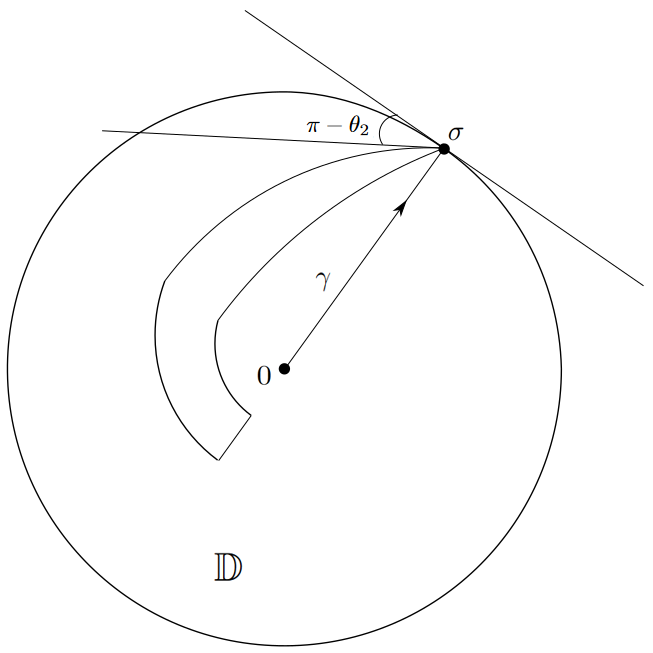}}\\
			\subfloat[\({[\theta_1,\theta_2]}={[\theta,\pi-\theta]}\)]{\includegraphics[scale=0.3]{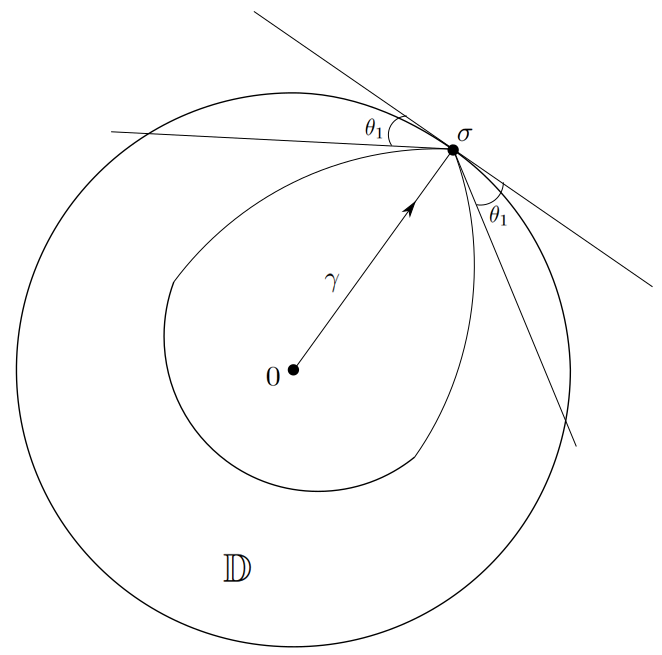}} 
			& \subfloat[\({[\theta_1,\theta_2]}={[\theta_1,\frac{\pi}{2}]}\)]{\includegraphics[scale=0.3]{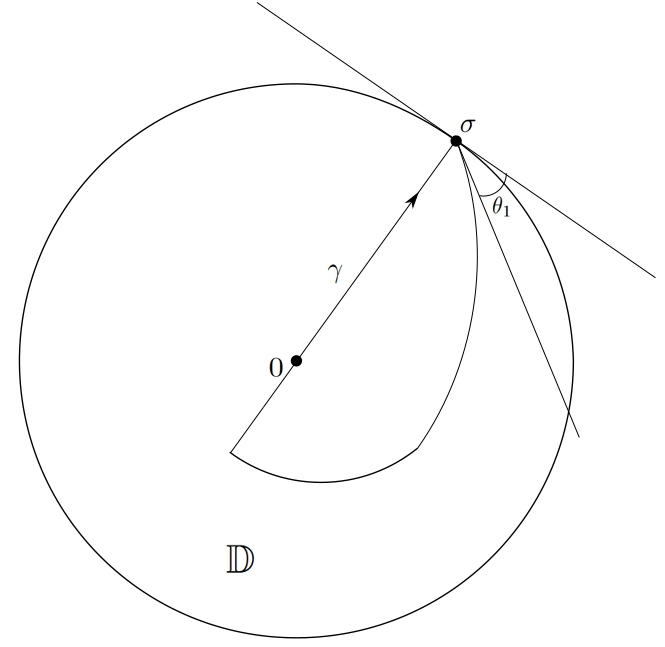}}\\
			\subfloat[\({[\theta_1,\theta_2]}={[\frac{\pi}{2},\theta_2]}\)]{\includegraphics[scale=0.36]{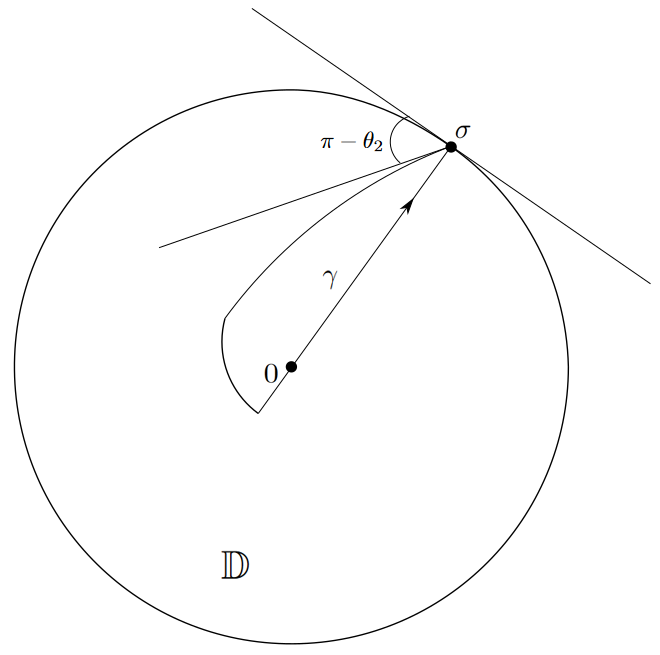}} 
			& \subfloat[\(\theta_1<\frac{\pi}{2}<\theta_2\)]{\includegraphics[scale=0.36]{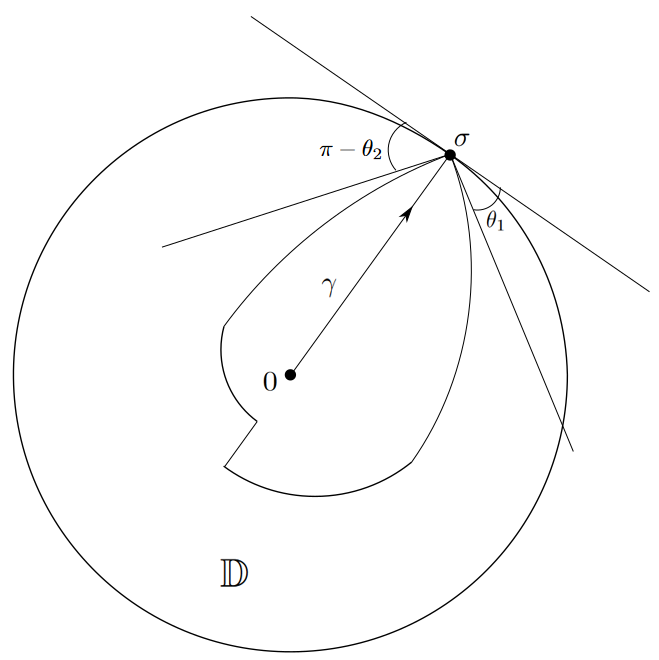}}\\
		\end{tabular}
		
	\end{tabularx}
	\caption{The set $A_\mathbb{D}(\gamma,\theta_1,\theta_2)$}
\end{figure}

\par In Figure 3, we have chosen the radius tending to $\sigma$ as the geodesic $\gamma$ and we present all the different versions of the set $A_\mathbb{D}(\gamma,\theta_1,\theta_2)$ depending on the angles $\theta_1,\theta_2$.

\par We now give the corresponding general definitions. Let $\Delta\subsetneq\mathbb{C}$ be a simply connected domain, $f:\mathbb{D}\to\Delta$ a Riemann map and $\gamma:[0,+\infty)\to\Delta$ a geodesic for the hyperbolic distance in $\Delta$. We can define the set $\tilde{A}_\Delta(\gamma,\theta_1,\theta_2)$ in a similar manner, using this time the hyperbolic distance in $\Delta$, corresponding hyperbolic sectors and the images $f(\mathbb{D}_\sigma^+),f(\mathbb{D}_\sigma^-)$, where this time, $\sigma$ is the point at which the preimage $f^{-1}\circ\gamma$ intersects the unit circle or equivalently $\lim\limits_{t\to+\infty}f^{-1}(\gamma(t))=\sigma$. For example, for $[\theta_1,\theta_2]\subset(0,\frac{\pi}{2})$, we have 
$$\tilde{A}_\Delta(\gamma,\theta_1,\theta_2)=\left[S_\Delta(\gamma,R(\frac{\pi}{2}-\theta_1))\setminus S_\Delta(\gamma,R(\frac{\pi}{2}-\theta_2))\right]\cap f(\mathbb{D}_\sigma^+),$$
and then $A_\Delta(\gamma,\theta_1,\theta_2)=\overline{\tilde{A}_\Delta(\gamma,\theta_1,\theta_2)}\cap\Delta$.

\begin{defin}
	Let $\Delta\subsetneq\mathbb{C}$ be a simply connected domain and $\{z_n\}$ a sequence in $\Delta$. We say that $\{z_n\}$ \textit{exhausts} the set $A_\Delta(\gamma,\theta_1,\theta_2)$, where $\gamma:[0,+\infty)\to\Delta$ is a geodesic for the hyperbolic distance in $\Delta$ and $\theta_1,\theta_2\in(0,\pi)$ if:
	\begin{enumerate}[(i)]
		\item  all the sets $A_\Delta(\gamma,\omega_1,\omega_2)$, where $\omega_1<\theta_1,\omega_2>\theta_2$, eventually contain $\{z_n\}$
		\item  for any $\theta\in[\theta_1,\theta_2]$ and any $\epsilon\in(0,\min\{\theta_1,\pi-\theta_2\})$, there is a subsequence of $\{z_n\}$ that is eventually contained in $A_\Delta(\gamma,\theta-\epsilon,\theta+\epsilon)$.
	\end{enumerate}
	This definition extends naturally to continuous curves as well.
\end{defin}

\par Next, we need the following proposition. During its proof, by $[z_1,z_2]$ we will mean the line segment joining the complex numbers $z_1,z_2$.
	\begin{prop}
	Let $\Delta\subsetneq\mathbb{C}$ be a simply connected domain and $f:\mathbb{D}\to\Delta$ a Riemann map. Suppose that $\mathbb{H}+a\subset\Delta\subseteq\mathbb{H}$, for some $a>0$. Let $\delta:[0,+\infty)\to\Delta$ be a continuous curve in $\Delta$ with $\lim\limits_{t\to+\infty}\delta(t)=\infty$. Suppose that $\delta$ exhausts $A_\mathbb{H}(\gamma,\theta_1,\theta_2)$, where $\gamma:[0,+\infty)\to\mathbb{H}$ is a geodesic for the hyperbolic distance in $\mathbb{H}$ such that $\gamma(0)=r_0>a$ and $\lim\limits_{t\to+\infty}\gamma(t)=\infty$. Then, there exists a $\sigma\in\partial\mathbb{D}$ such that $f^{-1}(\delta(t))$ converges by angle-set $[\theta_1,\theta_2]$ to $\sigma$, as $t\to+\infty$.
	\begin{proof}
		Let $f$ be a Riemann map for $\Delta$. By Carath\'{e}odory's Theorem, $f$ induces a homeomorphism $\tilde{f}:\mathbb{D}\cup\partial_C\mathbb{D}\to\Delta\cup\partial_C\Delta$. Also, there is a one-to-one correspondence between prime ends of the unit disk and points of the unit circle. Moreover, even though there might exist more than one prime ends of $\Delta$ corresponding to the point at infinity, there is a unique one that is determined by any chain of crosscuts that extends to the half-plane $\mathbb{H}+a$. Suppose that this unique prime end corresponds through $f^{-1}$ to a point $\sigma\in\partial\mathbb{D}$. Then, since $\delta([0,+\infty))$ is eventually contained in $\mathbb{H}+a$, it is necessary that $f^{-1}(\delta(t))\to\sigma$, as $t\to+\infty$. All that remains is to verify that the convergence is indeed by angle-set $[\theta_1,\theta_2]$.
		\par Since the geometry of $A_\mathbb{H}(\gamma,\theta_1,\theta_2)$ depends firmly on the selected angles, we will provide the proof in the case when $[\theta_1,\theta_2]\subset(0,\frac{\pi}{2})$. In the rest of the cases, the proof is almost the same, albeit with some minor modifications. Our hypotheses directly imply that the image of the geodesic $\gamma$ is the subset $[r_0,+\infty)$ of the real line. Using Lemma 3.1, we find that the boundary of $A_\mathbb{H}(\gamma,\theta_1,\theta_2)$ consists of  the two rays $\{re^{i(\theta_1-\frac{\pi}{2})}:r\ge r_0\},\;\{re^{i(\theta_2-\frac{\pi}{2})}:r\ge r_0\}$, the two lower halves $B_1,B_2$ of the circular arcs that join the points $r_0e^{i(\theta_1-\frac{\pi}{2})},r_0e^{i(\frac{\pi}{2}-\theta_1)}$ and $r_0e^{i(\theta_2-\frac{\pi}{2})},r_0e^{i(\frac{\pi}{2}-\theta_2)}$ respectively and a segment $[r_1,r_2]\subset\mathbb{R}$, with $0<r_1<r_2<r_0$, where $\{r_1\}=B_1\cap\mathbb{R}$ and $\{r_2\}=B_2\cap\mathbb{R}$.
		
		\begin{figure}
			\centering
			\includegraphics[scale=0.5]{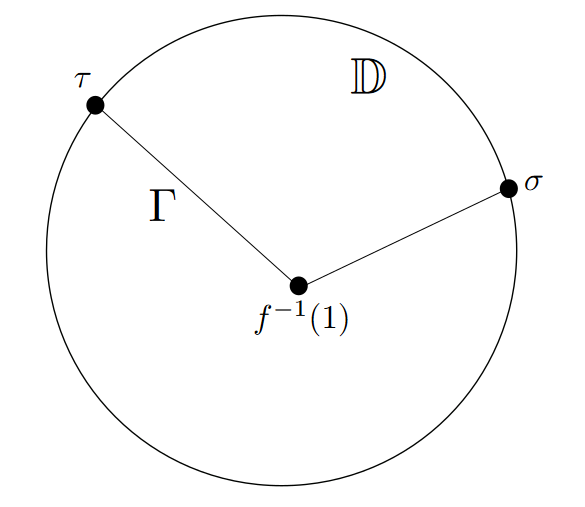}
			\caption{The set $\Gamma$}

		\end{figure}
		\par Without loss of generality, we can suppose that $a=1$ and that $\partial\Delta\cap\partial\mathbb{H}=\partial\Delta\cap(\partial\mathbb{H}+1)=\emptyset$ (in the case where $\Delta=\mathbb{H}$, the proof becomes much more straightforward).
		\par Set $x_0=\inf\{x\in\mathbb{R}\cap\partial\Delta\}$. Since $\partial\Delta\cap\partial\mathbb{H}=\emptyset$ and therefore $\Delta\subsetneq\mathbb{H}$, such a $x_0$ exists and, in fact, $x_0\in(0,1)$. Moreover, there exists a connected component $L$ of $\partial\Delta$ such that $x_0\in L, \inf\{\Im z:z\in L\}=-\infty$ and $\sup\{\Im z:z\in L\}=+\infty$. If either $\inf\{\Im z:z\in L\}\in(-\infty,0]$ or $\sup\{\Im z:z\in L\}\in[0,+\infty)$, then either $\inf\{x\in\mathbb{R}\cap\partial\Delta\}<x_0$ or $(-\infty,x_0)\subset\Delta$. Both cases lead to a contradiction and consequently there must exist such a component $L$. Therefore, there is a prime end $\xi$ of $\Delta$ corresponding to $x_0$ (if there are more than one, we just choose one randomly). It is possible that the impression $I[\xi]$ of this prime end is not a singleton and that $x_0$ is not accessible. In this case, since $I[\xi]\subset L$, we slightly change the construction and take as $x_0$ an accessible point belonging to $I[\xi]$ (for more on prime ends see \cite[Chapter 9]{prime}). As a result, through $f^{-1}$, the point $x_0\in\partial\Delta$ corresponds to a point $\tau\in\partial\mathbb{D}\setminus\{\sigma\}$.
		\par Consider in the unit disk $\mathbb{D}$ the set $\Gamma=(\tau,f^{-1}(1)]\cup[f^{-1}(1),\sigma)$ (see Figure 4). Obviously, $\Gamma$ separates the unit circle $\partial\mathbb{D}$ into two connected arcs. Passing to $\Delta$ through $f$, it is clear that $f(\Gamma)$ joins $x_0$ and $1$ and separates the prime ends of $\Delta$ into two sets $\partial_C\Delta^{+}$ and $\partial_C\Delta^-$. Denote as $\partial_C\Delta^+$ the subset of $\partial_C\Delta$ that contains the prime end of $\Delta$ determined by any chain of crosscuts with one end on the upper half of $L$ (i.e. the subset of $L$ with one end on $x_0$ and the other at infinity, but with imaginary parts tending to $+\infty$). Then, $\partial_C\Delta^-=\partial_C\Delta\setminus\partial_C\Delta^+$. We denote by $\partial\Delta^+$ and $\partial\Delta^-$ the sets of points of $\partial\Delta$ that correspond to the prime ends of $\partial_C\Delta^+$ and $\partial_C\Delta^-$ respectively. By our construction, $\partial\Delta^+$ and $\partial\Delta^-$ correspond through $f^{-1}$ to circular arcs of $\partial\mathbb{D}$ with common ends.
		
		\begin{figure}
			\centering
			\includegraphics[scale=0.6]{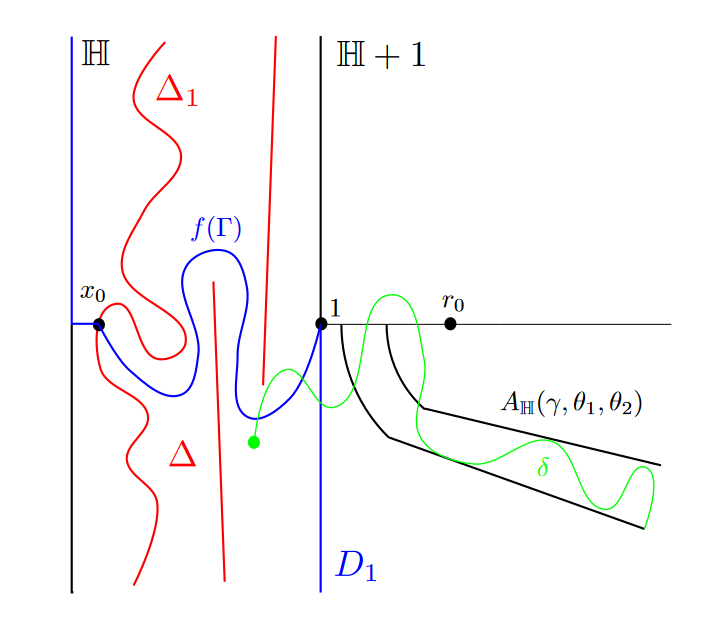}
			\caption{The sets of Proposition 3.1}
		\end{figure}

		\par For the rest of the proof we denote by $D_1$ the simply connected domain bounded by $\{z:\arg z=\frac{\pi}{2}\}\cup[0,x_0]\cup f((\tau,f^{-1}(1)])\cup[\{z:\arg z=-\frac{\pi}{2}\}+1]$ and by $D_2$ the simply connected domain bounded by $\{z:\arg z=-\frac{\pi}{2}\}\cup[0,x_0]\cup f((\tau,f^{-1}(1)])\cup[\{z:\arg z=\frac{\pi}{2}\}+1]$. Then, we set $\Delta_1:=D_1\cap\Delta$ and $\Delta_2:=D_2\cap\Delta$ (see Figure 5). By the construction above, both $\Delta_1$ and $\Delta_2$ are simply connected domains.
		\par Now, we will make use of the harmonic measure to prove that the convergence is indeed by angle-set $[\theta_1,\theta_2]$. By our hypothesis, there certainly exists a $t_0>0$, such that $\Re(\delta(t))>1$, for all $t>t_0$. Then, by the monotonicity property of the harmonic measure, we see that
		$$\omega(\delta(t),\partial\Delta^+,\Delta_1)\le\omega(\delta(t),\partial\Delta^+,\Delta),$$
		for all $t>t_0$. Moreover, by the maximum principle for harmonic functions,
		$$\omega(\delta(t),\partial\Delta^+,\Delta_1)\ge\omega(\delta(t),\{\arg z=\frac{\pi}{2}\},D_1),$$
		for all $t>t_0$. Combining the two and letting $t\to+\infty$, we find that
		\begin{equation}
			\limsup\limits_{t\to+\infty}\omega(\delta(t),\partial\Delta^+,\Delta)\ge\limsup\limits_{t\to+\infty}\omega(\delta(t),\{\arg z=\frac{\pi}{2}\},D_1).
		\end{equation}
	
		\par Then, we need to evaluate the quantity $\omega(\delta(t),\{\arg z=\frac{\pi}{2}\},D_1)$. Since $D_1\subset\mathbb{H}$, by the Strong Markov Property, for all $t>t_0$,
		\begin{eqnarray}
			\omega(\delta(t),\{\arg z=\frac{\pi}{2}\},\mathbb{H})&=&\omega(\delta(t),\{\arg z=\frac{\pi}{2}\},D_1)+ \\&+&\int\limits_{\partial D_1\setminus\partial\mathbb{H}}\omega(s,\{\arg z=\frac{\pi}{2}\},\mathbb{H})\cdot\omega(\delta(t),ds,D_1). \nonumber
		\end{eqnarray}

	\par Of course, $\partial D_1\setminus\partial\mathbb{H}=(0,x_0]\cup f((\tau,f^{-1}(1)])\cup[\{z:\arg z=-\frac{\pi}{2}\}+1]$. Consequently, the last integral can be broken into the two integrals
	
	$$I_1(t)=\int\limits_{(0,x_0]\cup f((\tau,f^{-1}(1)])}\omega(s,\{\arg z=\frac{\pi}{2}\},\mathbb{H})\cdot\omega(\delta(t),ds,D_1),$$
	
	$$I_2(t)=\int\limits_{[\{\arg z=-\frac{\pi}{2}\}+1]}\omega(s,\{\arg z=\frac{\pi}{2}\},\mathbb{H})\cdot\omega(\delta(t),ds,D_1).$$
	
	\par We remind that these integrals are valid for $t>t_0$. Since the harmonic measure has $1$ as an upper bound, it is directly computed that 
	$$I_1(t)\le\omega(\delta(t),(0,x_0]\cup f((\tau,f^{-1}(1)]),D_1).$$
	\par Without loss of generality, we can assume that there exists a line segment $A$ that joins $1$ with $\{z:\arg z=\frac{\pi}{2}\}$ such that $A\subset D_1$ (except, of course, for the endpoints). In case $f(\Gamma)$ reaches $1$ tangentially with respect to the vertical line $\{z:\Re z=1\}$ or even intersects $\mathbb{H}+1$, we can simply consider another $a>1$ such that $f((\tau,f^{-1}(a)])$ does not intersect $\mathbb{H}+a$ and is not tangent to the vertical line $\{z:\Re z=a\}$. This is possible, because $\tau\neq\sigma$ and the sets $f^{-1}(\mathbb{H}+a)$ are decreasing as $a\to+\infty$, with $\bigcap\limits_{a>1}\overline{f^{-1}(\mathbb{H}+a)}=\sigma$. Therefore, we can find a suitable $a>1$ such that, in the unit disk $\mathbb{D}$, the line segment $(\tau,f^{-1}(a)]$ does not intersect $f^{1}(\mathbb{H}+a)$ and is not tangent to $f^{1}(\mathbb{H}+a)$ at the point $f^{-1}(a)$. Then, we take $\Gamma=(\tau,f^{-1}(a)]\cup[f^{-1}(a),\sigma)$ and the construction of the line segment $A$ becomes well defined, while the proof can be executed for $\mathbb{H}+a\subset\Delta\subseteq\mathbb{H}$ without altering the result. 
	\par Denote by $H_1$ the simply connected domain bounded by $\{z:\arg z=\frac{\pi}{2}\},A$ and $[\{z:\arg z=-\frac{\pi}{2}\}+1]$ (see Figure 6). Obviously $H_1\subset D_1$. Then, by the maximum principle for harmonic functions, 
	$$\omega(\delta(t),(0,x_0]\cup f((\tau,f^{-1}(1)]),D_1)\le\omega(\delta(t),A,H_1),$$
	for all $t>t_0$. We can see that 
	\begin{eqnarray}
		\omega(\delta(t),A,H_1)&=&\omega(\delta(t)-1,A-1,H_1-1) \nonumber \\
		&\le&\omega(\delta(t)-1,A-1,\{-\frac{\pi}{2}<\arg z<\phi\}), \nonumber
	\end{eqnarray}
	where $\phi$ is the angle that $A$ forms with the real axis. Using consecutively the conformal mappings $g_1(z)=e^{i(\frac{\pi}{4}-\frac{\phi}{2})}z$, $g_2(z)=z^{\frac{2\pi}{\pi+2\phi}}$ and $g_3(z)=iz$, we map $\{z:-\frac{\pi}{2}<\arg z<\phi\}$ to the upper half-plane. It is easily checked that $g_3\circ g_2\circ g_1$ maps $A-1$ onto a line segment $[\alpha,\beta]$ on the negative semi-axis and $\delta(t)-1$ to $\delta_1(t)=i[(\delta(t)-1)e^{i(\frac{\pi}{4}-\frac{\phi}{2})}]^\frac{2\pi}{\pi+2\phi}$. Then, by the conformal invariance of the harmonic measure,
	$$\omega(\delta(t)-1,A-1,\{-\frac{\pi}{2}<\arg z<\phi\})=\omega(\delta_1(t),[\alpha,\beta],\{\Im z>0\}).$$
	 Combining all the above, we get
	$$\limsup\limits_{t\to+\infty}I_1(t)\le\limsup\limits_{t\to+\infty}\omega(\delta_1(t),[\alpha,\beta],\{\Im z>0\}).$$ 
	
		\begin{figure}
			\centering
			\includegraphics[scale=0.6]{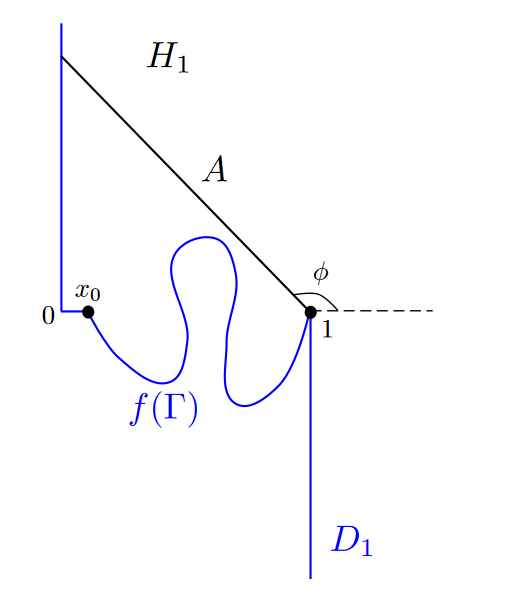}
			\caption{Construction of $A$ and $H_1$}
		\end{figure}
	Since $\delta_1(t)$ is eventually contained in an angle and converges to $\infty$, by our formula the last limit actually exists and is equal to 0, which leads to $\lim\limits_{t\to+\infty}I_1(t)=0$. 
	\par Next, we move to $I_2(t)$. Let $s\in[\{\arg z=-\frac{\pi}{2}\}+1]$ or $s=1-iy$, for some $y>0$. Again, by a known formula,
	$$\omega(s,\{\arg z=\frac{\pi}{2}\},\mathbb{H})=\frac{\arg(1-iy)+\frac{\pi}{2}}{\pi},$$
	with the last quantity tending to 0 as $y\to+\infty$. Let $\epsilon>0$. Because of the limit, there is a $y_0>0$ such that 
	$$\omega(1-iy,\{\arg z=\frac{\pi}{2}\},\mathbb{H})<\epsilon,$$
	for all $y>y_0$. Then, 
	\begin{eqnarray}
		I_2(t)&\le&\int\limits_{[1,1-iy_0]}\omega(s,\{\arg z=\frac{\pi}{2}\},\mathbb{H})\cdot\omega(\delta(t),ds,D_1)+\epsilon \nonumber \\
		&\le&\omega(\delta(t),[1,1-iy_0],D_1)+\epsilon. \nonumber
	\end{eqnarray}
	\par With the exact same technique as before, we use a conformal mapping onto the upper half-plane to show that
	$$\lim\limits_{t\to+\infty}\omega(\delta(t),[1,1-iy_0],D_1)=0.$$
	This means that $\limsup\limits_{t\to+\infty}I_2(t)\le\epsilon$, for all positive $\epsilon$, something that in turn leads to $\lim\limits_{t\to+\infty}I_2(t)=0$.
	Returning to (2) and remembering (1), we get 
	\begin{eqnarray}
		\limsup\limits_{t\to+\infty}\omega(\delta(t),\{\arg z=\frac{\pi}{2}\},\mathbb{H})&=&\limsup\limits_{t\to+\infty}\omega(\delta(t),\{\arg z=\frac{\pi}{2}\},D_1) \nonumber \\ &\le& \limsup\limits_{t\to+\infty}\omega(\delta(t),\partial\Delta^+,\Delta). \nonumber
	\end{eqnarray}
	
		\par The same is obviously true for the $\liminf$. We know from the given formulas that for every $\zeta$ on the ray $\{\arg z=\theta\}$, the harmonic measure $\omega(\zeta,\{\arg z=\frac{\pi}{2}\},\mathbb{H})$ remains constant and equal to $\frac{\theta+\frac{\pi}{2}}{\pi}$. Given that $\delta$ is a continuous curve that exhausts $A_\mathbb{H}(\gamma,\theta_1,\theta_2)$ and by the shape of $A_\mathbb{H}(\gamma,\theta_1,\theta_2)$, it follows that $\liminf\limits_{t\to+\infty}\arg(\delta(t))=\theta_1-\frac{\pi}{2}$ and $\limsup\limits_{t\to+\infty}\arg(\delta(t))=\theta_2-\frac{\pi}{2}$. As a result,
		\begin{eqnarray*}
	    \limsup\limits_{t\to+\infty}\omega(\delta(t),\partial\Delta^+,\Delta)&\ge&\limsup\limits_{t\to+\infty}\omega(\delta(t),\{\arg z=\frac{\pi}{2}\},\mathbb{H})\\&=&\limsup\limits_{t\to+\infty}\frac{\arg\delta(t)+\frac{\pi}{2}}{\pi}\\&=&\frac{\theta_2}{\pi}.
	    \end{eqnarray*}
    	Likewise,
	    $$\liminf\limits_{t\to+\infty}\omega(\delta(t),\partial\Delta^+,\Delta)\ge\frac{\theta_1}{\pi}.$$
	    Following a similar process, but this time with the use of $\{z:\arg z=-\frac{\pi}{2}\}$, $\partial\Delta^-$, $D_2$, $\Delta_2$, we find that
	    $$\limsup\limits_{t\to+\infty}\omega(\delta(t),\partial\Delta^-,\Delta)\ge\frac{\pi-\theta_1}{\pi},$$
	    $$\liminf\limits_{t\to+\infty}\omega(\delta(t),\partial\Delta^-,\Delta)\ge\frac{\pi-\theta_2}{\pi}.$$
	    However, since harmonic measure is always at most $1$, we have
	    $$\limsup\limits_{t\to+\infty}\omega(\delta(t),\partial\Delta^+,\Delta)+\liminf\limits_{t\to+\infty}\omega(\delta(t),\partial\Delta^-,\Delta)\le1=\frac{\theta_2}{\pi}+\frac{\pi-\theta_2}{\pi},$$
	    $$\liminf\limits_{t\to+\infty}\omega(\delta(t),\partial\Delta^+,\Delta)+\limsup\limits_{t\to+\infty}\omega(\delta(t),\partial\Delta^-,\Delta)\le1=\frac{\theta_1}{\pi}+\frac{\pi-\theta_1}{\pi}.$$
	    Consequently, all inequalities actually become equalities. Since $\delta$ is continuous, it is necessary that the cluster set of $\omega(\delta(t),\partial\Delta^+,\Delta)$ is equal to $[\frac{\theta_1}{\pi},\frac{\theta_2}{\pi}]$. Through $f^{-1}$, the set $\partial\Delta^+$ corresponds to a circular arc $B$ of $\partial\mathbb{D}$ with one end at $\sigma$ and the other at $\tau$. By conformal invariance, 
	    $$\liminf\limits_{t\to+\infty}\omega(f^{-1}(\delta(t)),B,\mathbb{D})=\frac{\theta_1}{\pi},$$
	    $$\limsup\limits_{t\to+\infty}\omega(f^{-1}(\delta(t)),B,\mathbb{D})=\frac{\theta_2}{\pi}.$$
	    We are now done, because in view of Remark 2.1, the last two equalities imply that as $t\to+\infty$, $f^{-1}(\delta(t))$ converges to $\sigma$ by angle-set $[\theta_1,\theta_2]$.
	\end{proof}
\end{prop}

\par With some minor changes, a similar proof can be executed for any horizontal geodesic of the right half-plane, since vertical translations are conformal automorphisms of $\mathbb{H}$ and so, they do not affect harmonic measure. We can think of it as a parallel move that does not affect the angles.

\begin{figure}
	\centering
	\includegraphics[scale=0.7]{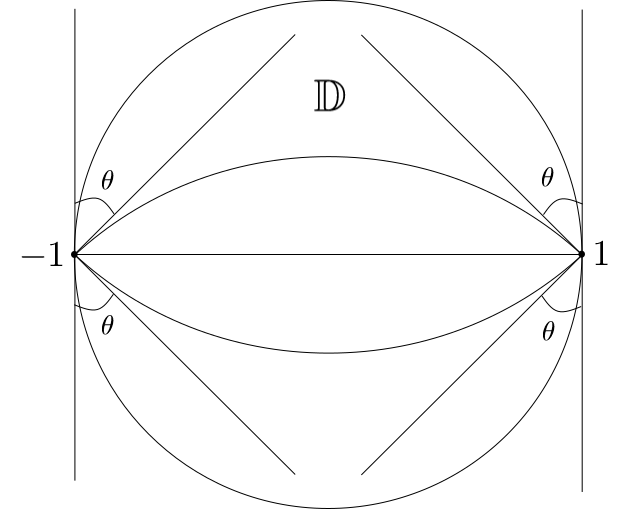}
	\caption{Hyperbolic sector in the unit disk}
\end{figure}

\par The last step, before proving the main theorem of this section, is to talk about hyperbolic sectors in the unit disk. Specifically, we will make some simple geometric considerations concerning hyperbolic sectors around a diameter. We can suppose that the diameter is the segment $(-1,1)$. Then a hyperbolic sector around this diameter (see Figure 7) is bounded by two circular arcs that join $-1$ and $1$ and intersect the unit circle with angles $\theta$ and $\pi-\theta$, where $\theta\in(0,\frac{\pi}{2})$. Restricting ourselves to a part of the diameter, for example the radius tending to 1, it is easy to see that the corresponding hyperbolic sector is bounded by two circular arcs $B_1,B_2$ that tend to $1$ and a third circular arc that joins $B_1$ and $B_2$ and is symmetric with respect to $(-1,1)$. Therefore, such a hyperbolic sector is eventually equivalent to a Stolz region and contains sequences that can converge to $1$ only by angle-set that is a subset of $(\theta,\pi-\theta)$. Keeping this in mind, we can understand how the sets $A_\mathbb{D}(\gamma,\theta_1,\theta_2)$ and the sequences that exhaust them behave.

\par We are now ready to prove our main result. For the sake of convenience, we will state the result again.

\begin{manualtheorem}{1.1}
	Let $\Delta\subsetneq\mathbb{C}$ be a simply connected domain and $f:\mathbb{D}\to\Delta$ a Riemann map. Let $\{z_n\}\subset\Delta$ be a sequence with no accumulation points in $\Delta$. Then, there exists a $\sigma\in\partial\mathbb{D}$ such that $\{f^{-1}(z_n)\}$ converges by angle-set $[\theta_1,\theta_2]\subset(0,\pi)$ to $\sigma$ if and only if there exist a simply connected domain $U\subsetneq\mathbb{C}$, a prime end $\xi\in\partial_C U$, $R>0$ and a geodesic $\gamma:[0,+\infty)\to U$ for the hyperbolic distance in $U$ with $\lim\limits_{t\to+\infty}\gamma(t)=\xi$ in the Carath\'{e}odory topology of $U$ such that
	\begin{enumerate}[\normalfont(i)]
		\item $E_U(\xi,R)\subset\Delta\subseteq U,$
		\item $\gamma(0)\in\Delta,$
		\item $\{z_n\}$ exhausts the set $A_U(\gamma,\theta_1,\theta_2).$
	\end{enumerate}
	\begin{proof}
	``$\implies$'' First, suppose that $\{f^{-1}(z_n)\}$ converges by angle-set $[\theta_1,\theta_2]$ to some $\sigma\in\partial\mathbb{D}$. By conformal invariance, we can take $\Delta=U=\mathbb{D},\xi=\sigma$ and $f=id_\mathbb{D}$. Therefore, condition (i) is obviously satisfied, for all $R>0$. Moreover, now $\{f^{-1}(z_n)\}\equiv\{z_n\}$, so $\{z_n\}$ converges by angle set $[\theta_1,\theta_2]$ to $\sigma$.	Consider $\gamma:[0,1)\to\mathbb{D}$ to be part of the radius that ends at $\sigma$. Condition (ii) is satisfied as well. Finally, by the way that the numbers $R_\mathbb{D}(\frac{\pi}{2}-\theta_i),\;i=1,2,$ were chosen for the definition of the set $A_\mathbb{D}(\gamma,\theta_1,\theta_2)$, the equivalence between hyperbolic sectors and Stolz regions and the geometry of the hyperbolic sectors around a radius, condition (iii) holds too.   \\
	``$\impliedby$'' Conversely, assume that conditions (i),(ii) and (iii) are satisfied. By conformal invariance, we can suppose that $U=\mathbb{H}$ and $\xi$ is the prime end of $\infty$ determined by the positive semi-axis. Then, there is an $a>0$ depending on $R$, such that $E_\mathbb{H}(\infty,R)=\mathbb{H}+a$. Therefore, condition (i) becomes $\mathbb{H}+a\subset\Delta\subseteq\mathbb{H}$. Also, the geodesic $\gamma$ of $\mathbb{H}$ satisfies $\lim\limits_{t\to+\infty}\gamma(t)=\infty$, which means that $\gamma$ must be a horizontal ray. Then, by the definition of exhaustion and Proposition 3.1, condition (iii) ensures that $\{f^{-1}(z_n)\}$ converges by angle-set $[\theta_1,\theta_2]$ to some $\sigma\in\partial\mathbb{D}$, as $t\to+\infty$.
	\end{proof}
\end{manualtheorem}

\begin{rem}
Examining the hyperbolic sectors around horizontal geodesics of $\mathbb{H}$, we can see that the existence of a suitable geodesic in the statement of Theorem 1.1 is equivalent to stating that the result is true for every such geodesic.
\end{rem}

\begin{rem}
	Let $\theta\in(0,\pi)$.
Suppose that $\{z_n\}$ is a sequence of points in the unit disk $\mathbb{D}$ with no accumulation points inside $\mathbb{D}$. Let $\sigma\in\partial\mathbb{D}$ and $\gamma:[0,+\infty)\to\mathbb{D}$ be the radius tending to $\sigma$ (or, even, any geodesic for the hyperbolic distance in the unit disk such that $\lim\limits_{t\to+\infty}\gamma(t)=\sigma$). We assume that $\{z_n\}$ is eventually contained in $A_\mathbb{D}(\gamma,\theta_1,\theta_2)$, for all $\theta_1\in(0,\theta)$ and all $\theta_2\in(\theta,\pi)$. Then, with the help of the conformal invariance of the harmonic measure and our formulas for the right half-plane, it is clear to see that $\{z_n\}$ necessarily converges to $\sigma$ and in fact, by angle $\theta$.
\par Conversely, if $\{z_n\}$ converges to $\sigma$ by angle $\theta$, then it is easily verified that $\{z_n\}$ is eventually contained in every set $A_\mathbb{D}(\gamma,\theta_1,\theta_2)$, for all $\theta_1\in(0,\theta)$ and all $\theta_2\in(\theta,\pi)$.
\end{rem}

Through a Riemann map, we can extend Remark 3.2 to any simply connected domain $\Delta\subsetneq\mathbb{C}$ in a suitable way. Therefore, a direct corollary of Theorem 1.1 is the following.

\begin{cor}
	Let $\Delta\subsetneq\mathbb{C}$ be a simply connected domain and $f:\mathbb{D}\to\Delta$ a Riemann map. Let $\{z_n\}\subset\Delta$ be a sequence with no accumulation points in $\Delta$. Then, there exists a $\sigma\in\partial\mathbb{D}$ such that $\{f^{-1}(z_n)\}$ converges by angle $\theta\in(0,\pi)$ to $\sigma$ if and only if there exist a simply connected domain $U\subsetneq\mathbb{C}$, a prime end $\xi\in\partial_C U, R>0$ and a geodesic $\gamma:[0,+\infty)\to U$ for the hyperbolic distance in $U$ with $\lim\limits_{t\to+\infty}\gamma(t)=\xi$ in the Carath\'{e}odory topology of $U$ such that
	\begin{enumerate}[\normalfont(i)]
		\item $E_U(\xi,R)\subset\Delta\subseteq U$,
		\item $\gamma(0)\in\Delta$,
		\item $\{z_n\}$ is eventually contained in the set $A_U(\gamma,\theta_1,\theta_2)$, for all $\theta_1\in(0,\theta)$ and all $\theta_2\in(\theta,\pi)$.
	\end{enumerate}
\end{cor}

\par Suppose, as above, that $\{z_n\}$ is eventually contained in the set $A_U(\gamma,\theta_1,\theta_2)$, for all $\theta_1\in(0,\frac{\pi}{2})$ and all $\theta_2\in(\frac{\pi}{2},\pi)$. Going back to the definition of $A_U(\gamma,\theta_1,\theta_2)$ and keeping in mind that $R(0)=0$, our assumption directly implies that $\lim\limits_{n\to+\infty}k_U(z_n,\gamma)=0$. On the other hand, one can see that if $\lim\limits_{n\to+\infty}k_U(z_n,\gamma)=0$, then $\{z_n\}$ is eventually contained in every hyperbolic sector around $\gamma$. Thus, it is also eventually contained in every set of the form $A_U(\gamma,\theta_1,\theta_2)$, where $0<\theta_1<\frac{\pi}{2}<\theta_2<\pi$. As a result, in the case that $\theta=\frac{\pi}{2}$ and the convergence is orthogonal, Corollary 3.1 is actually equivalent with Theorem 1.1 in \cite{orth} (for its statement see Section 1). Thus, Corollary 3.1 is a generalization for the convergence of a sequence by a certain angle, providing a necessary and sufficient condition.

\begin{rem}
	As we mentioned in the Introduction, the idea to tackle the problem of convergence by angle or angle-set in a simply connected domain $\Delta$ is the use of a larger simply connected domain $U$ that has a simpler geometry and is easier to work with. A natural question that arises after the proof of Theorem 1.1 and Corollary 3.1 is if the behavior of the sequence $\{z_n\}$ is the same with respect to both $\Delta$ and $U$. In particular, consider a Riemann map $f_\Delta:\mathbb{D}\to\Delta$ and a Riemann map $f_U:\mathbb{D}\to U$ and suppose that $E_U(\xi,R)\subset\Delta\subseteq U$, for some prime end $\xi\in\partial_C U$ and some $R>0$. Since $\Delta$ contains such a horodisk of $U$, it follows that $\xi\in\partial_C\Delta$ as well. Suppose that $\xi$ corresponds through $f_\Delta^{-1}$ to a point $\sigma_\Delta\in\partial\mathbb{D}$ and through $f_U^{-1}$ to a point $\sigma_U\in\partial\mathbb{D}$. Then, is it true that $\{f_\Delta^{-1}(z_n)\}$ converges by angle-set $[\theta_1,\theta_2]$ (or by angle $\theta$) to $\sigma_\Delta$ if and only if $\{f_U^{-1}(z_n)\}$ converges by angle-set $[\theta_1,\theta_2]$ (or by angle $\theta$) to $\sigma_U$?
	Actually, the answer to this question is positive. By conformal invariance, we can assume that $U=\mathbb{H}$ and that $\xi$ is the unique prime end of $\mathbb{H}$ corresponding to $\infty$. Then, following the proof of Proposition 3.1, we remember that in the end,
	$$\limsup\limits_{t\to+\infty}\omega(\delta(t),\{\arg z=\frac{\pi}{2}\},\mathbb{H})=\limsup\limits_{t\to+\infty}\omega(\delta(t),\partial\Delta^+,\Delta).$$
	The same is true for the liminf as well. Moreover, the positive semi-axis separates $\partial\mathbb{H}$ into the two connected components $\{\arg z=\frac{\pi}{2}\}$ and $\{\arg z=-\frac{\pi}{2}\}$, which through $f_\mathbb{H}^{-1}$ correspond to two circular arcs $B_\mathbb{H}^+$ and $B_\mathbb{H}^-$ respectively, that meet at $\sigma_\mathbb{H}$. Therefore, by the conformal invariance of the harmonic measure, we have
	$$\limsup\limits_{t\to+\infty}\omega(f_\mathbb{H}^{-1}(\delta(t)),B_\mathbb{H}^+,\mathbb{D})=\limsup\limits_{t\to+\infty}\omega(f_\Delta^{-1}(\delta(t)),B,\mathbb{D}).$$ 
	Again, the same is true for the liminf. As a result, in view of Remark 2.1, we get the desired result for the continuous curve $\delta$. Thus, the analogous result is true for the sequence $\{z_n\}$ as well.
\end{rem}

\section{Applications to Semigroups}

In this last part of the article, we connect the above work with the study of continuous semigroups of holomorphic functions of the unit disk (or from now on \textit{semigroups} in $\mathbb{D}$) and examine a consequence of Proposition 3.1. A semigroup in $\mathbb{D}$ is a family of holomorphic functions $\phi_t:\mathbb{D}\to\mathbb{D},t\ge0,$ with the properties:

\begin{enumerate}[(i)]
	\item $\phi_0$ is the identity in $\mathbb{D}$,
	\item $\phi_{t+s}=\phi_t\circ\phi_s$, for all $t,s\ge0$,
	\item $\lim\limits_{t\to s}\phi_t(z)=\phi_s(z)$, for all $s\ge0$ and all $z\in\mathbb{D}$.
\end{enumerate} 
\par An extensive presentation of the theory of semigroups can be found in the books \cite{abate},\cite{book},\cite{elin}.
\par When the functions of the semigroup $(\phi_t)$ do not have any fixed points, then $(\phi_t)$ is called \textit{non-elliptic}. For a non-elliptic semigroup $(\phi_t)$ in $\mathbb{D}$, there exists a unique point $\tau\in\partial\mathbb{D}$, namely the \textit{Denjoy-Wolff point} of the semigroup, such that for every $z\in\mathbb{D}$,
$$\lim\limits_{t\to+\infty}\phi_t(z)=\tau.$$
\par Usually, we assume without loss of generality, that the Denjoy-Wolff point is $1$.
\par An important tool connected to the study of a semigroup $(\phi_t)$ is the \textit{Koenigs function} $h$ of the semigroup. This function is the unique conformal mapping $h:\mathbb{D}\to\mathbb{C}$ satisfying $h(0)=0$ and
$$\phi_t(z)=h^{-1}(h(z)+t),$$ for all $t\ge0$ and all $z\in\mathbb{D}$. The simply connected domain $\Omega=h(\mathbb{D})$ is called the \textit{associated planar domain} of the semigroup and plays a great role in the classification of semigroups. It is convex in the positive direction (also known as starlike at infinity), namely $\Omega+t\subset\Omega$, for every $t\ge0$.
\par Going over to the classification of semigroups, a semigroup $(\phi_t)$ is called \textit{hyperbolic} if its associated planar domain $\Omega$ is contained in a horizontal strip. Otherwise, it is called \textit{parabolic}. In addition, we say that a parabolic semigroup is of \textit{zero hyperbolic step} if for some (equivalently for every) $s>0$ and some (equivalently every) $z\in\mathbb{D}$,
$$\lim\limits_{t\to+\infty}k_\mathbb{D}(\phi_t(z),\phi_{t+s}(z))=0.$$ If this limit is always positive, then we say that the parabolic semigroup is of \textit{positive hyperbolic step}. Consequently, we can divide non-elliptic semigroups in three classes: hyperbolic, parabolic of zero hyperbolic step and parabolic of positive hyperbolic step. For a study on the classification of semigroups in terms of the geometry of $\Omega$, we refer to \cite{class},\cite[Chapter 9]{book}. 
\par An interesting topic in the study of semigroups in $\mathbb{D}$ is the so-called ``slope'' problem. In this problem, we examine the trajectories of the semigroups (namely the sets $\{\phi_t(z):t\ge0\},z\in\mathbb{D}$) and most importantly the cluster set of $\arg(1-\bar{\tau}\phi_t(z))$, as $t\to+\infty$, where $\tau$ is the Denjoy-Wolff point of the semigroup. It has been proven (see \cite{hyper}) that in case of a hyperbolic semigroup, $\{\phi_t(z)\}$ always converges non-tangentially to the Denjoy-Wolff point, as $t\to+\infty$, while in case of a parabolic semigroup of positive hyperbolic step, the convergence is always tangential. In the third case, when the semigroup $(\phi_t)$ is parabolic of zero hyperbolic step, the situation is rather chaotic. At first, it was conjectured that the slope is a singleton, something that was proven under some additional assumptions in \cite{parab1,parab2}. However, in general, the convegence can be tangential or non-tangential and the cluster set of $\arg(1-\bar{\tau}\phi_t(z))$, as $t\to+\infty$, can be either a singleton $\{\theta\}$ or a set $[\theta_1,\theta_2]$, but it is always a compact connected subset of $[-\frac{\pi}{2},\frac{\pi}{2}]$. In \cite{bets,gume} examples are constructed where the slope is the whole interval $[-\frac{\pi}{2},\frac{\pi}{2}]$, while in \cite{non-tan,kelg} examples are constructed where the slope is a closed subinterval of $(-\frac{\pi}{2},\frac{\pi}{2})$. Furthermore, in \cite{orth}, the authors find geometric conditions that guarantee that the slope reduces to $\{0\}$. Finally, in \cite{asympt}, the authors characterize the type of convergence (tangential or non-tangential) of the trajectories of a semigroup $(\phi_t)$ to its Denjoy-Wolff point in terms of the shape of the associated planar domain.
\par In the present article, we provide some geometric conditions on the associated planar domain $\Omega$ that are sufficient for the slope to be a singleton $\{\frac{\pi}{2}-\theta\},\;\theta\in(0,\pi)$. For this purpose, we will first prove the next proposition which is very similar to Proposition 3.1, but concerns the convergence by a certain angle and not by an angle-set. In particular, this next proposition is a generalization of Proposition 4.1 in \cite{orth}. We set $U_\theta=\{z:-\theta<\arg z<\pi-\theta\}$, where $\theta\in(0,\pi)$. Actually, for $\theta=\frac{\pi}{2}$, the statements of the two propositions are identical. However, the proofs are entirely different. The proof in \cite{orth} is based on hyperbolic geometry, while this next proof relies on harmonic measure. We can consider Proposition 4.1 as a subcase of Proposition 3.1. In fact, the following proof is directly implied by the one of Proposition 3.1.
\begin{prop}
	Let $\Delta\subsetneq\mathbb{C}$ be a simply connected domain and $f:\mathbb{D}\to\Delta$ a Riemann map. We suppose that $U_\theta+a\subset\Delta\subseteq U_\theta$, for some $a>0$. Then, there exists a $\sigma\in\partial\mathbb{D}$ such that $f^{-1}(t)$ converges by angle $\theta$ to $\sigma$, as $t\to+\infty$.
	\begin{proof}
		
		Using the rotation $g(z)=e^{i(-\frac{\pi}{2}+\theta)}z$ which is a conformal mapping and, as such, keeps the harmonic measure invariant, we can assume that $\mathbb{H}+a\subset\Delta\subseteq\mathbb{H}$. Consider the continuous curve $\delta:[0,+\infty)\to\Delta$ with $\delta(t)=e^{i(-\frac{\pi}{2}+\theta)}(t+a)$, which is actually a reparametrization and then a rotation through $g$ of the positive semi-axis. Obviously, $\lim\limits_{t\to+\infty}\delta(t)=\infty$. By the geometry of the right half-plane and its hyperbolic sectors and by the definition of the sets $A_\mathbb{H}(\gamma,\theta_1,\theta_2)$, it is easy to check that $\delta([0,+\infty))$ is eventually contained in every set of the form $A_\mathbb{H}(\gamma,\theta_1,\theta_2)$, where $\gamma$ is the geodesic for the hyperbolic distance of the right half-plane that has as its image the positive semi-axis and $0<\theta_1<\theta<\theta_2<\pi$. Therefore, with the exact same proof as in Proposition 3.1, we find that $\delta(t)$ converges to $\infty$ by angle $\theta$, as $t\to+\infty$. Rotating back through $g^{-1}$ and returning to the unit disk through $f^{-1}$, we get the desired result.

	\end{proof}
\end{prop}

\par In the following corollary, $\Omega$ will not be contained in any horizontal strip and so the semigroup cannot be hyperbolic. In addition, we already said that the slope will not reduce to $\{-\frac{\pi}{2}\}$ or $\{\frac{\pi}{2}\}$ or in other words the convergence will not be tangential. Therefore, the semigroup will not be parabolic of positive hyperbolic step. Having restricted ourselves to non-elliptic semigroups, we can only talk about a parabolic semigroup $(\phi_t)$ of zero hypebolic step. For $\alpha_1,\alpha_2\in(0,\pi]$, consider the sector 
$$U(\alpha_1,\alpha_2)=\{z:-\alpha_1<\arg z<\alpha_2\}.$$

\begin{cor}
	Let $(\phi_t)$ be a parabolic semigroup in $\mathbb{D}$ of zero hyperbolic step, $h$ the corresponding Koenigs function, $\Omega=h(\mathbb{D})$ and $\tau\in\partial\mathbb{D}$ the Denjoy-Wolff point of $(\phi_t)$. If there exist $\alpha_1,\alpha_2\in(0,\pi]$ with $\alpha_1+\alpha_2\ge\pi$ and $a>0$ such that $U(\alpha_1,\alpha_2)+a\subset\Omega\subseteq U(\alpha_1,\alpha_2)$, then $\lim\limits_{t\to+\infty}\arg(1-\bar{\tau}\phi_t(z))=\frac{\pi}{2}\cdot\frac{\alpha_2-\alpha_1}{\alpha_1+\alpha_2}$, for every $z\in\mathbb{D}$.
	\begin{proof}
	Let $z\in\mathbb{D}$ such that $h(z)\in\mathbb{R}$. Surely, the trajectory $\phi_t(z)$ converges to the Denjoy-Wolff point $\tau$, as $t\to+\infty$. Through the Koenigs function $h$, all the trajectories correspond to horizontal rays that converge to $\infty$ in the positive direction. So, we have a curve $h(\phi_t(z))\subset\mathbb{R},t\ge0,$ with $\lim\limits_{t\to+\infty}h(\phi_t(z))=\infty$ and $U(\alpha_1,\alpha_2)+a\subset\Omega\subseteq U(\alpha_1,\alpha_2)$. Applying consecutively the conformal mappings $g_1(z)=e^{i\frac{\alpha_1-\alpha_2}{2}}z$, $g_2(z)=z^{\frac{\pi}{\alpha_1+\alpha_2}}$ and $g_3(z)=e^{-i\frac{\pi}{2}\cdot\frac{\alpha_1-\alpha_2}{\alpha_1+\alpha_2}}z$, we can map $U(\alpha_1,\alpha_2)$ conformally onto the set $U=\{z:-\frac{\pi\alpha_1}{\alpha_1+\alpha_2}<\arg z<\frac{\pi\alpha_2}{\alpha_1+\alpha_2}\}$, while also $g_t(z)=g_3\circ g_2\circ g_1\circ h\circ\phi_t(z)\in\mathbb{R}$, for all $t\ge0$. We can easily see that the set $U$ is actually a half-plane and as a result we can directly make use of Proposition 4.1 to find that $g_t(z)$ converges by angle $\frac{\pi\alpha_1}{\alpha_1+\alpha_2}$ to $\infty$, as $t\to+\infty$. Therefore, keeping in mind the conformal invariance of the harmonic measure, $\phi_t(z)$ converges by angle $\frac{\pi\alpha_1}{\alpha_1+\alpha_2}$ to $\tau$, as $t\to+\infty$. By definition, this means that $\lim\limits_{t\to+\infty}\arg(1-\bar{\tau}\phi_t(z))=\frac{\pi}{2}\cdot\frac{\alpha_2-\alpha_1}{\alpha_1+\alpha_2}$. Finally, we know (see e.g. \cite[p.492]{book}) that the dynamic behavior of every trajectory is the same. As a consequence, the desired result holds for every $z\in\mathbb{D}$.
	\end{proof}
\end{cor}
\begin{rem}
	In the statement of this last corollary, we demand that $\alpha_1+\alpha_2\ge\pi$. This added condition is needed because in the case when $\alpha_1+\alpha_2<\pi$, the chain of conformal mappings we used, does not guarantee that Proposition 4.1 can be applied.
\end{rem}
It is worth mentioning that a more general result was proven in \cite{slope}, where the authors use functions which measure the angular displacement of the boundary of the associated planar domain with respect to a fixed vertical straight line. In fact, the above Corollary follows immediately from Corollary 5.11 in \cite{slope}.

\section*{Acknowledgements}
I thank Professor D. Betsakos, my advisor, for his advice during the preparation of this work. I also thank the referee for his/her very valuable suggestions.
\par This research did not receive any specific grant from funding agencies in the public, commercial, or not-for-profit sectors.
\\
Declarations of interest: none.

\end{document}